\documentclass[9pt,onecolumn]{article}
\usepackage{definitions}

\title{First-Order Methods for Linear Programming}
\author{Haihao Lu}
\date{}

\begin{document}

\maketitle

% \begin{abstract}
%     Linear programming is one of the most
% \end{abstract}

\section{Introduction}
Linear programming (LP)~\cite{kantarovich1939mathematical, dantzig1948programming, dantzig1963linear,vanderbei2020linear,schrijver1998theory,luenberger1984linear} is the seminal optimization problem that has spawned and grown into today’s rich and diverse optimization modeling and algorithmic landscape. LP is used in just about every arena of the global economy, including transportation, telecommunications, production and operations scheduling, as well as in support of strategic decision-making~\cite{hazell1974competitive,delson1992linear,dahleh1994control,liu2008choice,charnes1959application,zhou2008linear}. Eugene Lawler is quoted as stating in 1980 that LP “is used to allocate resources, plan production, schedule workers, plan investment portfolios and formulate marketing (and military) strategies. The versatility and economic impact of linear optimization models in today’s industrial world is truly awesome.”

% It refers to an optimization problem that minimizes/maximizes an affine objective function over a polytope, defined by a set of linear equality and inequality constraints~\cite{kantarovich1939mathematical, dantzig1948programming, dantzig1963linear,vanderbei2020linear,schrijver1998theory,luenberger1984linear}. It is arguably one of the most important classes of optimization problems with prevalent applications in every aspect of industrial civilization, including agriculture, engineering, revenue management, finance, and the energy industry~\cite{hazell1974competitive,delson1992linear,dahleh1994control,liu2008choice,charnes1959application,zhou2008linear}, just to name a few. 

% \cite{charnes1954stepping,hanssmann1960linear,bowman1956production,manne1960linear,liu2008choice,anderson2000hotel}
 % The field of mathematical programming perhaps started in late 1940s when 

% Over the years, LP has become one of the most fundamental tools in applied math, operation research and computer science, finding wide-ranging applications in diverse fields including allocation, scheduling, inventory control, chip design, network flow, revenue management, and many others~\cite{charnes1954stepping,hanssmann1960linear,bowman1956production,manne1960linear,anderson2000hotel,liu2008choice}. 

Since the late 1940s, the exceptional modeling capabilities of LP have spurred extensive investigations into efficient algorithms for solving LP. Presently, the most widely recognized methods for LP are Dantzig's simplex method~\cite{dantzig1963linear,dantzig1990origins} and interior-point methods~\cite{karmarkar1984new,renegar1988polynomial,monteiro1989interior,wright1997primal}. The state-of-the-art commercial LP solvers, which are based on variants of these two methods, are quite mature, and can reliably deliver extremely accurate solutions. However, it is extremely challenging to further scale either of these two algorithms beyond the problem sizes they can currently handle. More specifically, the computational bottlenecks of both methods involve matrix factorization to solve linear equations, which leads to two fundamental challenges as the size of the problem increases:
\begin{itemize}
    \item In the commercial LP solvers, the simplex method utilizes LU-factorization and the interior point method  utilizes Cholesky factorization when solving the linear equations. While the constraint matrix is extremely sparse in practice, it is often the case that the factorization is much denser than the constraint matrix. This is the reason that commercial LP solvers require more memory than just storing the LP instance itself, and it may lead to out-of-memory errors when solving large instances. Even fitting in memory, the factorization may take a long time for large instances.
    \item It is highly challenging to take advantage of modern computing architectures, such as GPUs or distributed systems, to solve the linear equations in these two methods. For simplex, multiple cores are generally employed only in the ``pricing'' step that selects variables entering or leaving the basis, and the speedups are generally minimal after two or three cores~\cite{huangfu2018parallelizing}. For interior point methods, commercial solvers can use multiple threads to factor the associated linear systems, and speedups are generally negligible after at most six shared-memory cores, although further scaling is occasionally possible.
    \end{itemize}
Another classic approach to solving large-scale LPs is to use decomposition algorithms, such as Dantzig-Wolfe decomposition, Benders decomposition, etc. While these algorithms are useful for solving some problems, they are tied to certain block structures of the problem instances and suffer from slow tail convergence, and thus they are not suitable for general-purpose LP solvers.

Given the above limitations of existing methods, first-order methods (FOMs) have become increasingly attractive for large LPs. FOMs only utilize gradient information to update their iterates, and the computational bottleneck is matrix-vector multiplication, as opposed to matrix factorization. Therefore, one only needs to store the LP instance in memory when using FOMs rather than storing any factorization. Furthermore, thanks in part to the recent development of machine learning/deep learning, FOMs scale very well on GPU and distributed computing platforms.

Indeed, using FOMs for LP is not a new idea. As early as the 1950s, initial efforts were made to employ FOMs for solving LP problems. Notably, with the intuition to make big jumps rather than ``crawling along edges", \cite{brown1951computational} discussed steepest ascent to maximize the linear objective under linear inequality constraints. \cite{zoutendijk1960methods,zoutendijk1970some} pioneered the development of feasible direction methods for LP, where the iterates consistently move along a feasible descent direction. The subsequent development of steepest descent gravitational methods in~\cite{chang1989steepest} also falls within this category. Another early approach to first-order methods for solving LP is the utilization of projected gradient algorithms~\cite{rosen1961gradient,lemke1961constrained}. All these methods, however, still require solving linear systems to determine the appropriate direction for advancement. Solving the linear systems that arise during the update can be highly challenging for large instances. Moreover, the computation of projections onto polyhedral constraints can be arduous and even intractable for large-scale instances.

% \subsection{Overview of recent FOMs for LP}
The recent surge of interest in large-scale applications of LP has prompted the development of new FOM-based LP algorithms, aiming to further scale up/speed up LP. The four main solvers are:
% \begin{itemize}
    % \item 
    
    {\bf PDLP~\cite{applegate2021practical}.} PDLP utilizes primal-dual hybrid gradient method (PDHG) as its base algorithm and introduces practical algorithmic enhancements, such as presolving, preconditioning, adaptive restart, adaptive choice of step size, and primal weight, on top of PDHG. Right now, it has three implementations:  a prototype implemented in Julia (\href{https://github.com/google-research/FirstOrderLp.jl}{FirstOrderLp.jl}) for research purposes, a production-level C++ implementation that is included in Google \href{https://developers.google.com/optimization}{OR-Tools}, and an internal distributed version at Google. The internal distributed version of PDLP has been used to solve real-world problems with as many as 92B non-zeros~\cite{blog}, which is one of the largest LP instances that are solved by a general-purpose LP solver.
    
    % A distributed version of PDLP has demonstrated its efficacy in solving real-world LP instances with as many as 92 billion non-zero elements in the contraint matrix. This scale is significantly larger, by a factor of hundreds to thousands, than what state-of-the-art commercial LP solvers are capable of handling. 
    % \item 
    {\bf ABIP~\cite{lin2021admm,deng2022new}.} \href{https://github.com/leavesgrp/ABIP}{ABIP} is an ADMM-based IPM. The core algorithm of ABIP is still a homogeneous self-dual embedded interior-point method. Instead of approximately minimizing the log-barrier penalty function with a Newton's step, ABIP utilizes multiple steps of alternating direction method of multipliers (ADMM). The  $\mathcal O\left(\frac{1}{\epsilon}\log\left(\frac{1}{\epsilon}\right)\right)$ sublinear complexity of ABIP was presented in \cite{lin2021admm}. Recently, \cite{deng2022new} includes new enhancements, i.e., preconditioning, restart, hybrid parameter tuning, on top of ABIP (the enhanced version is called ABIP+). ABIP+ is numerically comparable to the Julia implementation of PDLP. ABIP+ now also supports a more general conic setting when the proximal problem associated with the log-barrier in ABIP can be efficiently computed.
    % \item 

    {\bf ECLIPSE~\cite{basu2020eclipse}.} ECLIPSE is a distributed LP solver designed specifically for addressing large-scale LPs encountered in web applications. These LPs have a certain decomposition structure, and the effective constraints are usually much less than the number of variables. ECLIPSE looks at a certain dual formulation of the problem, then utilizes accelerated gradient descent to solve the smoothed dual problem with Nesterov's smoothing. This approach is shown to have $\mathcal O(\frac{1}{\epsilon})$ complexity, and it is used to solve web applications with $10^{12}$ decision variables~\cite{basu2020eclipse} and real-world web applications at LinkedIn platform~\cite{ramanath2022efficient,acharya2023promoting}.
    
    % that employs accelerated gradient descent on smoothed dual LP. It is designed and motivated specifically for addressing large-scale LPs encountered in web applications with billions to trillions of decision variables and constraints. \cite{basu2020eclipse} provides proof of the sublinear convergence of ECLIPSE and demonstrates its effectiveness in solving problems involving $10^{12}$ decision variables.
    % \item 
    {\bf SCS~\cite{o2016conic,o2021operator}.}     
    Splitting conic solver \href{https://github.com/cvxgrp/scs}{(SCS)} tackles the homogeneous self-dual embedding of general conic programming using ADMM. As a special case of conic programming, SCS can also be used to solve LP. Each iteration of SCS involves projecting onto the cone and solving a system of linear equations with similar forms so that it only needs to store one factorization in memory. Furthermore, SCS supports solving the linear equations with an iterative method, which only uses matrix-vector multiplication. 
% \end{itemize}

In the rest of this article, we focus on discussing the theoretical and computational results of PDLP to provide a solid introduction to the field of FOMs for LP. Indeed, many theoretical guarantees we mention below can be applied directly to SCS, since ADMM is a pre-conditioned version of PDHG. Furthermore, many enhancements proposed in PDLP, such as preconditioning and restart, have been implemented in ABIP+.

% Many of the enhancements s

% theoretical guarantees  can also be applied to other solvers .

\section{PDHG for LP}
% PDHG is the base algorithm used in PDLP. In this section, we present the theoretical results of PDHG for LP. 

For this section, we look at the standard form of LP for the sake of simplicity, and most of these results can be extended to other forms of LP. More formally, we consider
% \begin{equation}\label{eq:lp-pd}
%     \begin{aligned}[c]
%     \min_{x\in \mathbb R^n}~~ &~ c^\top x \\
% \text{subject to:}~~ & ~ Ax = b,\; x \geq 0
%     \end{aligned}
%     \qquad\qquad\qquad
%     \begin{aligned}[c]
% \max_{y\in \mathbb R^{m}} \quad & b^T y \\
% \text{subject to:}\quad & A^Ty\leq c\ ,
%     \end{aligned}
% \end{equation} 
\begin{equation}\label{eq:primal}
    \begin{aligned}[c]
    \min_{x\in \mathbb R^n}~~ &~ c^\top x \\
\text{s.t.}~~ &~ Ax=b \\
& ~ x\geq 0 \ ,
    \end{aligned}
\end{equation}
% \begin{align}\label{eq:primal}
%     \begin{split}
%         & \min_{x}\;c^Tx \\
%         & \ \mathrm{s.t.}\; Ax=b \\
%         & \quad\quad x\geq 0 \ ,
%     \end{split}
% \end{align}
where $A\in\mathbb R^{m\times n}$, $b\in \mathbb R^m$ and $c\in \mathbb R^n$. The dual of standard form LP is given by
% \begin{align}\label{eq:dual}
%     \begin{split}
%         & \max_{y}\;b^Ty \\
%         & \ \mathrm{s.t.}\; A^Ty\leq c \ . 
%     \end{split}
% \end{align}

\begin{equation}\label{eq:dual}
    \begin{aligned}[c]
    \max_{y\in \mathbb R^m}~~ &~ b^\top y \\
\text{s.t.}~~ &~ A^\top y\leq c \ .
    \end{aligned}
\end{equation}

We start with discussing the most natural FOMs to solve LP as well as their limitations, which then lead to PDHG. Next, we discuss the convergence guarantees of vanilla PDHG for LP, and built upon that, we presented restarted PDHG, which can be shown as an optimal FOM for solving LP, i.e., it matches with the complexity lower bound. Then, we discuss how PDHG can detect infeasibility without additional effort.
\subsection{Initial attempts}

\begin{figure}[ht!]
	\centering\includegraphics[scale=0.4]{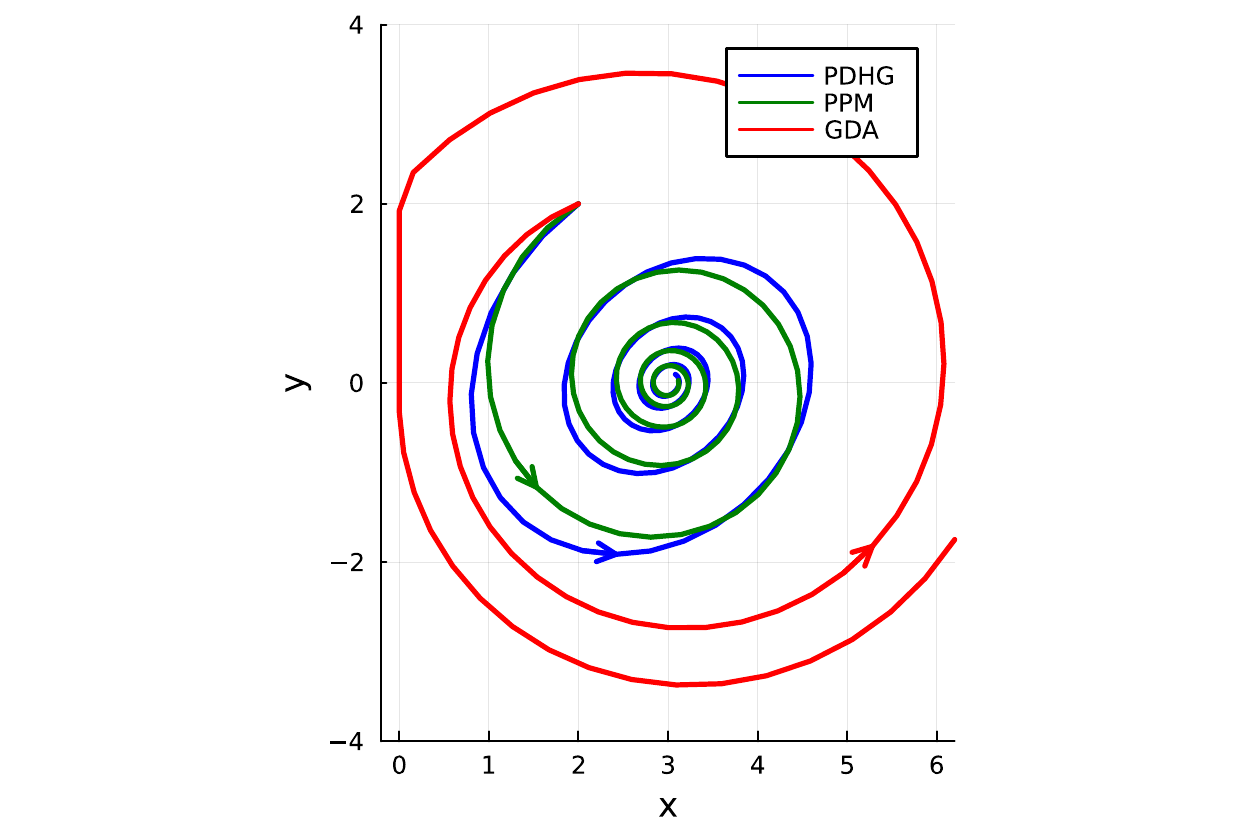}
	\caption{Trajectories of GDA, PPM and PDHG to solve a simple bilinear problem \eqref{eq:bilinear} with initial solution $(2, 2)$ and step-size $\eta=\sigma=0.2$.}
	\label{fig:pdhg-gda}
\end{figure}

The most natural FOM for solving a constrained optimization problem, such as LP \eqref{eq:primal}, is perhaps the projected gradient descent (PGD), which has iterated update:
\begin{equation*}\label{alg:pgd}
    x^{k+1}=\text{proj}_{\{x\in\mathbb R_+^n|Ax=b\}}(x^k-\eta c)\ ,
\end{equation*}
where $\eta$ is the step-size of the algorithm.
All the nice theoretical guarantees of PGD for general constrained convex optimization problems can be directly applied to LP. Unfortunately, computing the projection onto the constrained set (i.e., the intersection of an affine subspace and the positive orthant) involves solving a quadratic programming problem, which can be as hard as solving the original LP, and thus PGD is not a practical algorithm. To disentangle linear constraints and (simple) nonnegativity of variables, a natural idea is to dualize the linear constraints $Ax=b$ and consider the primal-dual form of the problem \eqref{eq:minmax}:
\begin{equation}\label{eq:minmax}
    \min_{x\geq 0}\max_y\; L(x,y):=c^Tx-y^TAx+b^Ty \ .
\end{equation}
Convex duality theory~\cite{boyd2004convex} shows that the saddle points to \eqref{eq:minmax} can recover the optimal solutions to the primal problem~\eqref{eq:primal} and the dual problem~\eqref{eq:dual}. 
For the primal-dual formulation of LP~\eqref{eq:minmax}, the most natural FOM  is perhaps the projected gradient descent-ascent method (GDA), which has iterated update:
\begin{equation*}\label{alg:gda}
    \begin{cases}
        x^{k+1}=\text{proj}_{\mathbb R^n_+}(x^k+\eta A^Ty^k-\eta c) \\ y^{k+1}= y^k-\sigma Ax^k+\sigma b\ ,
    \end{cases}
\end{equation*}
where $\eta$ and $\sigma$ are the primal and the dual step-size, respectively.
The projection of GDA is onto the positive orthant for the primal variables and it is cheap to implement. Unfortunately, GDA does not converge to a saddle point of \eqref{eq:minmax}. For instance, Figure \ref{fig:pdhg-gda} plots the trajectory of GDA on a simple primal-dual form of LP
\begin{equation}\label{eq:bilinear}
    \min_{x\ge 0}\max_y\; (x-3)y \ ,
\end{equation}
where $(3,0)$ is the unique saddle point. As we can see, the GDA iterates diverge and spins farther away from the saddle point. Thus, GDA is not a good algorithm for \eqref{eq:minmax}.

% We run gradient descent and ascent on the primal and dual variable respectively.  However, the problem hereby is that the iterates of GDA diverge even on an two-dimensional unconstrained bilinear problem. To illustrate, the red line in Figure \ref{fig:pdhg-gda} plots the trajectory of GDA to solve
% \begin{equation}\label{eq:bilinear}
%     \min_{x\ge 0}\max_y\; (x-3)y \ .
% \end{equation}
% The algorithm is initiated at $(0.5,0.5)$ and $(0,0)$ is the unique optimal solution. We can see the spiral pattern of GDA iterates but it spins farther away from the origin. Thus the direct extension of PGD to the primal-dual problems, i.e., GDA, does not work.

Another candidate algorithm to solve such primal-dual problem is the proximal point method (PPM), proposed in the seminal work of  Rockafellar~\cite{rockafellar1976monotone}.  PPM has the following iterated update:
\begin{equation}\label{eq:ppm}
    (x^{k+1},y^{k+1})\leftarrow \arg\min_{x\geq 0}\max_y\;L(x,y)+\frac{1}{2\eta}\|x-x^k\|_2^2-\frac{1}{2\sigma}\|y-y^k\|_2^2\ .
\end{equation}
Unlike GDA, PPM exhibits nice theoretical properties for solving primal-dual problems (see Figure \ref{fig:pdhg-gda} for example). However, its update rule is implicit and requires solving the subproblems arising in \eqref{eq:ppm}. This drawback makes PPM more of a conceptual rather than a practical algorithm.

To overcome these issues of GDA and PPM, we consider primal-dual hybrid gradient method (PDHG, a.k.a Chambolle-Pock algorithm)~\cite{chambolle2011first,zhu2008efficient}. 
PDHG is a first-order method for convex-concave primal-dual problems originally motivated by applications in image processing. In the case of LP, the update rule is straightforward:
% which has the following iterate updades:
% Primal-Dual Hybrid Gradient (PDHG, a.k.a Chambolle-Pock algorithm)~\cite{chambolle2011first,zhu2008efficient} is a primal-dual algorithm motivated by problems in image processing and later finds its applications beyond. PDHG solves \eqref{eq:minmax} with the following updates
\begin{equation}\label{alg:pdhg}
    \begin{cases}
        x^{k+1}\leftarrow \text{proj}_{\mathbb R^n_+}(x^k+\eta A^Ty^k-\eta c) \\ y^{k+1}\leftarrow y^k-\sigma A(2x^{k+1}-x^k)+\sigma b\ ,
    \end{cases}
\end{equation}
where $\eta$ is the primal step-size and $\sigma$ is the dual step-size. Similar to GDA, the algorithm alternates with the primal and the dual variables, and the difference is in the dual update, one utilizes the gradient at the extrapolated point $2x^{k+1}-x^k$. The extrapolation helps with the convergence of the algorithm, as we can see in Figure \ref{fig:pdhg-gda}. Indeed, one can show the PDHG is a preconditioned version of PPM (see the next section for more details), and thus share the nice convergence properties with PPM, but it does not require solving the implicit update~\ref{eq:ppm}. The computational bottleneck of PDHG is the matrix-vector multiplication (i.e., in $A^T y$ and $Ax$).
 
 % Unlike PPM, the computational bottleneck of PDHG is matrix-vector multiplication (i.e., in $A^T y$ and $Ax$). As we can see in Figure \ref{fig:pdhg-gda}, PDHG converges to the unique stationary point of the primal-dual problem~\eqref{eq:bilinear}. Indeed, one can show good properties of the convergence of PDHG for LP, which we present in the next section.

 % The primal update of PDHG is identical to GDA, but PDHG calculates the dual gradient at the extrapolated point $2x^{k+1}-x^k$ in contrast with at last iterates $x^k$ for GDA. It turns out that this extrapolation boosts the convergence of PDHG (see next section). As a sanity check, we also run PDHG on the simple unconstrained bilinear problem \eqref{eq:bilinear}. The blue line in Figure \ref{fig:pdhg-gda} refers to the iterates of PDHG on \eqref{eq:bilinear} which also have spiral pattern and converges rapidly to optimality.

% This puzzling convergent/divergent behaviors of PDHG and GDA can be explained via the difference in their high-resolution ODE~\cite{lu2022sr}.

\subsection{Theory of PDHG on LP}
In this section, we start with presenting the sublinear convergence rate for the average iterates and the last iterates of PDHG. Then we present the sharpness of the primal-dual formulation of LP that leads to the linear convergence of PDHG on LP.

% Then, we introduce a sharpness definition for primal-dual problems, and show that PDHG has linear convergence under this sharpness condition, which is satisfied by LP. 

For simplicity of the exposition, we assume the primal and dual step-sizes are equal, i.e., $\eta=\sigma=s$ throughout this section and the next section. This can be done without the loss of generality by rescaling the primal and the dual variables. Furthermore, we assume the LP instances are feasible and bounded in this and the next section, and we will discuss infeasibility detection later on. For notational simplicity, we denote $z=(x,y)$ as the pair of the primal and dual solution,
    $P_s=\begin{pmatrix}
        \frac{1}{s}I & A^T \\ A & \frac{1}{s}I
    \end{pmatrix}, $
and $\|z\|_{P_s}=\sqrt{\langle z,P_sz\rangle}$. 

The $P_s$ norm is the inherent norm for PDHG, and it plays a critical role in the theoretical analysis of PDHG. One can clearly see this by noticing the PDHG update~\eqref{alg:pdhg} can be rewritten as
\begin{equation}\label{eq:pdhg-new}
    P_s(z^k-z^{k+1})\in \mF(z^{k+1}) \ ,
\end{equation}
where $\mF(z)=\begin{pmatrix} \partial_x L(x,y) \\ -\partial_y L(x,y) \end{pmatrix}$ is the sub-differential of the objective. This also showcases that PDHG is a pre-conditioned version of PPM with $P_s$ norm by noticing that one can rewrite the update rule of PPM as $\frac{1}{s}(z^k-z^{k+1})\in \mF(z^{k+1})$. Armed with this understanding, one can easily obtain the sublinear convergence of of PDHG:

\begin{thm}[Average iterate convergence of PDHG {\cite{chambolle2011first,chambolle2016ergodic,lu2023unified}}]\label{thm:average}
    Consider the iterates $\{z^k\}_{k=0}^{\infty}$ of PDHG \eqref{alg:pdhg} for solving \eqref{eq:minmax}. Denote $\bar z^k=(\bar x^k,\bar y^k)=\frac{1}{k}\sum_{i=1}^k z^i$ as the average iterates. Then it holds for any $0<s\le \frac{1}{\|A\|_2}$, $k\ge 1$, and $z=(x,y)\in \mathcal Z$ that
    \begin{equation*}
        L(\bar x^k,y)-L(x,\bar y^k)\leq \frac{1}{2k}\|z-z^0\|^2_{P_s} \ .
    \end{equation*}
    %     where $P_s=\begin{pmatrix}
    %     \frac{1}{s}I & A^T \\ A & \frac{1}{s}I
    % \end{pmatrix}$.
\end{thm}

% We define $\|z\|_{P_s}=\langle z,P_sz\rangle$ as the canonical norm of PDHG with step-size $s$. This norm plays a critical role in the theory for PDHG. Note that $\mathrm{dist}(z,Z):=\min_{u\in Z}\|z-u\|_2$ and $\mathrm{dist}_{P_s}(z,Z):=\min_{u\in Z}\|z-u\|_{P_s}$.

% Proposition \ref{prop:prop-sublinear} shows the $\mathcal O(1/\sqrt{k})$ rate of the last iterates of PDHG on solving LP. 
\begin{thm}[Last iterate convergence of PDHG \cite{lu2022infimal}]\label{prop:prop-sublinear}
Consider the iterates $\{z_k\}_{k=0}^{\infty}$ of PDHG \eqref{alg:pdhg} for solving \eqref{eq:minmax}. Let $z^*\in \mathcal Z^*$ be an optimal point to \eqref{eq:minmax}. Suppose the step-size satisfies $s<\frac{1}{\Vert A \Vert}$. Then, it holds for any iteration $k\ge 1$, $z=(x,y)\in \mathcal Z$ and any optimal solution $z^*$ that
\begin{equation*}
    L(x^k,y)-L(x,y^k) \leq \frac{1}{\sqrt k}\pran{\Vert z^0-z^* \Vert_{P_s}^2+\Vert z^{0}-z^* \Vert_{P_s}\Vert z^*-z \Vert_{P_s}} \ .
\end{equation*}
\end{thm}

Theorem \ref{thm:average} and Theorem \ref{prop:prop-sublinear} show that the average iterates of PDHG have $\mathcal O(1/k)$ convergence rate, and the last iterates of PDHG have $\mathcal O(1/\sqrt k)$ convergence rate. The above two convergence results are not limited to LP and PDHG. \cite{lu2022infimal,lu2023unified} show that these results work for convex-concave primal-dual problems (i.e., for general $L(x,y)$ as long as $L$ is convex in $x$ and concave in $y$) and for more general algorithms as long as the update rule can be written as an instance of \eqref{eq:pdhg-new}, such as PPM, alternating direction method of multipliers (ADMM), etc.

% the $\mathcal O(1/k)$ rate of the average iterates of PDHG on solving LP.Theorem \ref{thm:average} and Theorem \ref{prop:prop-sublinear}
Theorem \ref{thm:average} and Theorem \ref{prop:prop-sublinear} also imply that the average iterates of PDHG have a faster convergence rate than the last iterates. Numerically, one often observes that the last iterates of PDHG exhibit faster convergence (even linear convergence) than the average iterates. This is due to the structure of LP, which satisfies a certain regularity condition that we call the sharpness condition. To formally define this condition, we first introduce a new progress metric, normalized duality gap, defined in \cite{applegate2022faster}:
% Next, we further show the linear convergence of PDHG. First we introduce a new progress metric, which we dub normalized duality gap, and discuss the sharpness of LP \eqref{eq:minmax}. Theorem \ref{thm:linear} presents the linear rate of the last iterates under sharpness.
\begin{mydef}[Normalized duality gap {\cite{applegate2022faster}}]
    For a primal-dual problem~\eqref{eq:minmax} and a solution $z=(x,y)\in \mathcal Z$, the normalized duality gap with radius $r$ is defined as
    \begin{equation}\label{eq:ndg}
    \rho_r(z)=\max_{\hat z \in W_r(z)}\frac{L(x,\hat y)-L(\hat x,y)}{r} \ ,
\end{equation}
where $W_r(z)=\{ \hat z\in\mathcal Z \;|\; \Vert z-\hat z\Vert_{2}\leq r \}$ is a ball centered at $z$ with radius $r$ intersected with $\mathcal Z=\mathbb R^n_+ \times \mathbb R^m$.
\end{mydef}
Normalized duality gap is a valid progress measurement for LP, since $\rho_r(z)$ is a continuous function and $\rho_r(z)=0$ if and only if $z$ is an optimal solution to \eqref{eq:minmax}. One can show that LP is indeed a sharp problem: the normalized duality gap $\rho_r(z)$ of \eqref{eq:minmax} is sharp in the standard sense.
\begin{prop}[{\cite{applegate2022faster}}]\label{prop:sharp}
    The primal-dual formulation of linear programming \eqref{eq:minmax} is $\alpha$-sharp on the set $\|z\|_2\leq R$ for all $r\leq R$, i.e., there exists a constant $\alpha>0$ and it holds for any $z$ with $\|z\|_2\leq R$ and any $r\leq R$ that 
    \begin{equation*}
        \alpha\mathrm{dist}(z,\mathcal Z^*)\leq \rho_r(z) \ ,
    \end{equation*}
    where $\mZ^*$ is the optimal solution set, $\mathrm{dist}(z,\mZ^*)=\min_{z^*\in\mZ^*}\|z-z^*\|$ is the distance between $z$ and $\mZ^*$.%, and $\mathrm{diam}(S)$ is the diameter of the set $S$.
\end{prop}

The next theorem shows that the last iterates of PDHG have global linear convergence on LP (more generally, sharp primal-dual problems). In contrast, the average iterates always have sublinear convergence.
\begin{thm}[{\cite{lu2022infimal}}]\label{thm:linear}
    Consider the iterates $\{z^k\}_{k=0}^{\infty}$ of PDHG \eqref{alg:pdhg} to solve \eqref{eq:minmax}. Suppose the step-size $s\leq\frac{1}{2\Vert A \Vert}$, and \eqref{eq:minmax} is $\alpha$-sharp on $B_R(0)$, where $R$ is the upper bound of $\|z^k\|_2$. Then, it holds for any iteration $k\ge \left\lceil 4e/(s\alpha)^2\right\rceil$ that
    \begin{equation*}
    dist_{P_s}(z^k,\mathcal Z^*)\leq \exp\pran{\frac 12-\frac{k}{2\left\lceil 4e/(s\alpha)^2\right\rceil}}dist_{P_s}(z^0,\mathcal Z^*) \ .
    \end{equation*}
\end{thm}
% \begin{rem}
%     Theorem \ref{thm:linear} implies that we need in total $\mathcal O\pran{\pran{\frac{1}{s\alpha}}^2\log\pran{\frac{1}{\epsilon}}}$ iterations of PDHG to find an $\epsilon$-close solution $z$.
% \end{rem}

\subsection{Optimal FOM for LP}\label{sec:restart}

Theorem \ref{thm:linear} shows that the last iterates of PDHG have linear convergence with complexity $\mathcal O\pran{\pran{\frac{\|A\|_2}{\alpha}}^2\log\pran{\frac{1}{\epsilon}}}$ when the optimal step-size is chosen. A natural question is whether there exists FOM with faster convergence for LP. The answer is yes, and it turns out a simple variant of PDHG achieves faster linear convergence than PDHG and it matches with the complexity lower bound. 

Algorithm \ref{alg:pdhg-restart} formally presents this algorithm, which we dub restarted PDHG. This is a two-loop algorithm. The inner loop runs PDHG until one of the restart conditions holds. At the end of each inner loop, the algorithm restarts the next outer loop at the running average of the current epoch.

% We here present an adaptive scheme: essentially, we restart 

% ere we propose an adaptive restart scheme that does not need an estimate of $\alpha$. We restart when the normalized duality gap has sufficient decay, i.e., we restart the algorithm if
%     \begin{equation}\label{eq:adaptive-restart}
%     \begin{cases}
%         \rho_{\|\bar z^{n,k}-z^{n,0}\|_{2}}(\bar z^{n,k}) \leq e^{-1} \rho_{\|z^{n,0}-z^{n-1,0}\|_{2}}(\bar z^{n,0}),\; n\geq 1 \\ 
%         k \geq \tau^0,\; n=0
%     \end{cases}
% \end{equation}
% where $\tau^0$ is the first restart interval length and can be selected by the user.

% In the previous section, we show the sublinear convergence of average iterates and linear convergence of last iterates under sharpness of \eqref{eq:minmax} for vanilla PDHG. In this section, we present a restart variant of PDHG (Algorithm \ref{alg:pdhg-restart}) and show restarted PDHG exhibits faster linear rate than vanilla PDHG. Two restart schemes are proposed: fixed frequency restart and adaptive restart. Moreover, it turns out that restarted PDHG is indeed the optimal first-order methods for solving LP. We prove this by presenting the lower bound result.

\begin{algorithm}%[H]
\caption{Restarted PDHG for \eqref{eq:minmax}}
\label{alg:pdhg-restart}
\SetKwInOut{Input}{Input}
\Input{Initial point $(x^0,y^0)$, step-sizes $0<s<\frac{1}{\|A\|_2}$, outer loop counter $n\leftarrow 0$.}

\Repeat{\upshape $(x^{n+1,0},y^{n+1,0})$ convergence}{
  initialize the inner loop counter $k\leftarrow0$;\\
  \Repeat{\upshape restart condition holds}{
    $x^{n,k+1}\leftarrow\text{proj}_{\mathbb R^n_+}(x^{n,k}+\eta A^Ty^{n,k}-\eta c)$;\\
    $y^{n,k+1}\leftarrow y^{n,k}-\sigma A(2x^{n,k+1}-x^{n,k})+\sigma b$;\\
    $(\bar x^{n,k+1}, \bar y^{n,k+1})\leftarrow\frac{1}{k+1}\sum_{i=1}^{k+1} (x^{n,i},y^{n,i})$;
  }
  initialize the initial solution $(x^{n+1,0}, y^{n+1,0})\leftarrow(\bar x^{n,k+1}, \bar y^{n,k+1})$;\\
  $n\leftarrow n+1$;
}
\end{algorithm}
A crucial component of the algorithm is when to restart. Suppose we know the sharpness constant $\alpha$ and $\|A\|_2$. The fixed frequency restart scheme proposed in \cite{applegate2022faster} is to restart the algorithm every
\begin{equation}\label{eq:fixed-restart}
    k^*=\left\lceil \frac{4e\|A\|_2}{\alpha} \right\rceil \ 
\end{equation}
iterations. The next theorem presents the linear convergence rate of PDHG with the above fixed frequency restart.

\begin{thm}[{\cite{applegate2022faster}}]\label{thm:linear-adaptive}
    Consider the iterates $\{z^{n,0}\}_{n=0}^{\infty}$ generated by Algorithm \ref{alg:pdhg-restart} with fixed frequency restart, namely, we restart the outer loop if \eqref{eq:fixed-restart} holds, for solving \eqref{eq:minmax}. Then for any $\epsilon>0$, the algorithm finds $z^{n,0}$ such that $\mathrm{dist}_{P_s}(z^{n,0},\mathcal Z^*)\le \epsilon$ within
    $$
    \mathcal O\pran{\frac{\|A\|_2}{\alpha}\log\pran{\frac{1}{\epsilon}}}
    $$
   PDHG iterations.
\end{thm}

Furthermore, \cite{applegate2022faster} shows that restarted PDHG matches with the complexity lower bound of a wide range of FOMs for LP. In particular, we consider the span-respecting FOMs:

% The rest of the section presents the lower bound of first-order methods for solving linear programming \eqref{eq:minmax}. Combined with Theorem \ref{thm:linear-fixed}, it shows that restarted PDHG is optimal for LP.
\begin{mydef}\label{def:span-respecting}
    An algorithm is span-respecting for an unconstrained primal-dual problem $\min_x\max_y\; L(x,y)$ if
    \begin{align*}
        \begin{split}
            & \ x^k\in x^0+\mathrm{span}\{\nabla_x L(x^i,y^j): \forall i,j\in\{1,...,k-1\}\}\\
            & \ y^k\in y^0+\mathrm{span}\{\nabla_y L(x^i,y^j): \forall i\in\{1,...,k\},\forall j\in\{1,...,k-1\} \} \ .
        \end{split}
    \end{align*}
\end{mydef}
Definition \ref{def:span-respecting} is an extension of the span-respecting FOMs for minimization~\cite{nesterov2003introductory} in the primal-dual setting. 
% If $L$ is bilinear then with appropriate indexing of their iterates, primal-dual algorithms including GDA, PDHG, extragradient~\cite{korpelevich1976extragradient} and their restarted variants satisfies Definition \ref{def:span-respecting}. 
Theorem \ref{thm:lb} provides a lower complexity bound of span-respecting primal-dual algorithms for LP.

\begin{thm}[Lower complexity bound {\cite{applegate2022faster}}]\label{thm:lb}
    Consider any iteration $k\geq 0$ and parameter value $\gamma>\alpha>0$. There exists an $\alpha$-sharp linear programming with $\|A\|_2=\gamma$ such that the iterates $z^k$ of any span-respecting algorithm  satisfies that
    \begin{equation*}
        \mathrm{dist}(z^k,\mathcal Z^*)\geq \pran{1-\frac{\alpha}{\gamma}}^{k} \mathrm{dist}(z^0,\mathcal Z^*) \ .
    \end{equation*}
\end{thm}

Together with Theorem \ref{thm:linear-adaptive}, Theorem \ref{thm:lb} shows that restarted PDHG is an optimal FOM for LP. 

In practice, one might not know the sharpness constant $\alpha$ and the smoothness constant $\gamma$. \cite{applegate2022faster} proposes an adaptive restart scheme, which essentially restarts the algorithm whenever the normalized duality gap has a constant factor shrinkage. The adaptive restart scheme does not require knowing the parameters of the problem, and it leads to a nearly optimal complexity (up to a $\log$ term).

% When to restart the algorithm is a key component of the algorithm. \cite{applegate2022faster} proposes two schemes: restart with fixed frequency \eqref{eq:fixed-restart} and adaptive restart (see \eqref{eq:adaptive-restart} in Section \ref{sec:enhancement}). In this section, we present the theoretical results of fixed frequency restart for the sake of simplicity. We would also like to mention that PDHG with adaptive restart also exhibits (nearly) optimal rate upto the log term~\cite{applegate2022faster}. 
%The extra log term can be eliminated using fixed restart frequency and hence provide an optimal convergence rate \cite{applegate2022faster}.
% \begin{rem}
%     Theorem \ref{thm:lb} implies the following lower complexity bound for any first-order method to achieve an $\epsilon$-accuracy solution: $$\Omega\pran{\frac{L}{\alpha}\log{\frac{1}{\epsilon}}} \ .$$
%     A direct consequence of Theorem \ref{thm:linear-fixed} and Theorem \ref{thm:lb} is the optimality of restarted PDHG among FOM for solving LP.
% \end{rem}

\subsection{Infeasibility detection}

The convergence results of PDHG in the previous section require the LP to be feasible and bounded. In practice, it is occasionally the case that LP is infeasible or unbounded, thus infeasibility detection is a necessary feature for any LP solver. In this section, we investigate the behavior of PDHG on infeasible/unbounded LPs, and claim that the PDHG iterates encode infeasibility information automatically. 

The easiest way to describe the infeasiblity detection property of PDHG is perhaps to look at it from an operator perspective. More formally, we use $T$ to represent the operator for one step of PDHG iteration, i.e., $z^{k+1}=T(z^k)$ where $T$ is specified by \eqref{alg:pdhg}. Next, we introduce the infimal displacement vector of the operator $T$, which plays a central role in the infeasibility detection of PDHG:
\begin{mydef}
    For the operator $T$ induced by PDHG on \eqref{eq:minmax},  we call $$v:=\arg\min_{z\in\text{range}(T-I)} \|z\|_2^2$$ its infimal displacement vector (which is uniquely defined~\cite{pazy1971asymptotic}).
\end{mydef}
It turns out that if LP is primal (or dual) infeasible, then the dual (or primal) variables diverge like a ray with direction $v$. Furthermore, the corresponding dual (or primal) part of $v$ provides an infeasibility certificate for the primal. Table \ref{tab:infeas} summaries such results. More formally, 

\begin{thm}[Behaviors of PDHG for infeasible LP {\cite{applegate2021infeasibility}}]\label{thm:full-infeas}
    Consider the primal problem \eqref{eq:primal} and dual problem \eqref{eq:dual}. Assume $s< \frac{1}{\|A\|}$, let $T$ be the operator induced by PDHG on \eqref{eq:minmax}, and let $\{z^k=(x^k,y^k)\}_{k=0}^{\infty}$ be a sequence generated by the fixed-point iteration from an arbitrary starting point $z^0$. Then, one of the following holds:

    (a). If both primal and dual are feasible, then the iterates $(x^k,y^k)$ converge to a primal-dual solution $z^*=(x^*,y^*)$ and $v=(T-I)(z^*)=0$.

    (b). If both primal and dual are infeasible, then both primal and dual iterates diverge to infinity. Moreover, the primal and dual components of the infimal displacement vector $v=(v_x,v_y)$ give certificates of dual and primal infeasibility, respectively.

    (c). If the primal is infeasible and the dual is feasible, then the dual iterates diverge to infinity, while the primal iterates converge to a vector $x^*$. The dual-component $v_y$ is a certificate of primal infeasibility. Furthermore, there exists a vector $y^*$ such that $v=(T-I)(x^*,y^*)$.

    (d). If the primal is feasible and the dual is infeasible, then the same conclusions as in the previous item hold by swapping primal with dual.
\end{thm}

\begin{table}[h]
\centering
{\large
\begin{tabular}{|c|c|c|}
\hline
\diagbox[width=10em]{\textbf{Primal}}{\textbf{Dual}}                   & \textbf{Feasible} & \textbf{Infeasible} \\ \hline
\textbf{Feasible}   &  $x^k,y^k$ both converge              &  $x^k$ diverges, $y^k$ converges    \\ \hline
\textbf{Infeasible} &  $x^k$ converges, $y^k$ diverges                 &  $x^k,y^k$ both diverge                    \\ \hline
\end{tabular}
}
\caption{Behavior of PDHG for solving \eqref{eq:minmax} under different feasibility assumptions.}
\label{tab:infeas}
\end{table}

Furthermore, one can show that the difference of iterates and the normalized iterates converge to the infimal displacement vector $v$ with sublinear rate:
\begin{thm}[{\cite{applegate2021infeasibility,davis2016convergence}}]\label{thm:infeas}
    Let $T$ be the operator induced by PDHG on \eqref{eq:minmax}. Then there exists a finite $z^*$ such that $T(z^*)=z^*+v$ and for any such $z^*$ and all $k$:

    (a) (Difference of iterates)
    \begin{equation*}
        \min_{j\leq k}\|v-(z^{j+1}-z^j)\|\leq \frac{1}{\sqrt k}\|z^0-z^*\|_{P_s}\ ,
    \end{equation*}

    (b) (Normalized iterates)
    \begin{equation*}
        \left\|v-\frac{1}{k}(z^k-z^0)\right\|_{P_s}\leq \frac{2}{k}\|z^0-z^*\|_{P_s}\ .
    \end{equation*}
\end{thm}

Theorem \ref{thm:full-infeas} and Theorem \ref{thm:infeas} show that the difference of iterates and the normalized iterates of PDHG can recover the infeasibility certificates with sublinear rate. While the normalized iterates have faster sub-linear convergence than the difference of iterates, one can show that the difference of iterates converge linearly to the infimal displacement vector under additional regularity conditions~\cite{applegate2021infeasibility}. In practice, PDLP periodically checks whether the difference of iterates or the normalized iterates provide an infeasibility certificate, and the performance of these two sequences is instance-dependent.

\section{PDLP}
In the previous section, we present theoretical results of PDHG for LP. In the solver PDLP, there are additional algorithmic enhancements on top of PDHG to boost the practical performance. In this section, we summarize the enhancements as well as the numerical performance of the algorithms. These results were based on the Julia implementation and were presented in~\cite{applegate2021practical}. The algorithm in the c++ implementation of PDLP is almost identical to the Julia implementation, with two minor differences: it supports two-sided constraints, and utilizes Glop presolve instead of Papilo presolve.

% the relevant theoretical results of PDLP. However, there is also a gap between the theoretical analysis and the practical solver. In this section, we present several algorithmic enhancements of PDLP to boost the practical performance. The experimental testing and comparison with other solvers are then provided~\cite{applegate2021practical}.

PDLP solves a more general form of LP:
\begin{equation}\label{eq:lp}
    \begin{aligned}[c]
    \min_{x\in \mathbb R^n}~~ &~ c^\top x \\
\text{s.t.}~~ &~ Gx \geq h \\
& ~ Ax = b \\
& ~ l \leq x \leq u
    \end{aligned}
    % \qquad\qquad\qquad
%     \begin{aligned}[c]
% \max_{y\in \mathbb R^{m_1 + m_2}, \lambda\in\mathbb R^n} \quad & q^T y + l^T \lambda^{+} - u^T \lambda^{-} \\
% \text{subject to:}\quad & c - K^T y = \lambda \\
% &  y_{1:m_1} \geq 0 \\
% & \lambda \in \Lambda\ ,
%     \end{aligned}
\end{equation} 
%with $n, m_1, m_2 \in \N$.
where $G \in \mathbb R^{m_1\times n}$, $A \in \mathbb R^{m_2 \times n}$, $c \in \mathbb R^{n}$, $h \in \mathbb R^{m_1}$, $b \in \mathbb R^{m_2}$, $l \in (\mathbb R \cup \{ -\infty \})^{n}$, $u \in (\mathbb R \cup \{ \infty \})^{n}$.
% $K := \begin{pmatrix}
% G \\
% A
% \end{pmatrix}$, and 
% $$
% \Lambda = \Lambda_{1} \times \dots \times \Lambda_{n} \quad \Lambda_i :=
% \begin{cases}
% \{ 0 \} & l_i = -\infty, ~ u_i = \infty,  \\
% \mathbb R^{-} & l_i = -\infty, ~ u_i \in \mathbb R \\
% \mathbb R^{+} & l_i \in \mathbb R, ~ u_i = \infty \\
% \mathbb R & \text{otherwise} 
% \end{cases}
% $$
% is the set of variables $\lambda$ such that the dual objective is finite. 
The primal-dual form of the problem becomes:
\begin{flalign}\label{eq:primal-dual}
\min_{x \in X} \max_{y \in Y} L(x,y) := c^T x - y^T K x + q^T y
\end{flalign}
where $K^\top = \begin{pmatrix} G^T, A^T \end{pmatrix}$ and $q^T := \begin{pmatrix}
h^T,
b^T
\end{pmatrix}$, $X := \{x \in \mathbb R^n : l \leq x \leq u \}$, and $Y := \{y \in \mathbb R^{m_1+m_2} : y_{1:m_1} \geq 0\}.$ In PDLP, the primal and the dual are reparameterized as
\begin{equation*}
    \eta = s/\omega,\; \sigma=s\omega\quad \text{with}\; s,\omega>0\ ,
\end{equation*}
where $s$ (which we call step-size) controls the scale of the step-size, and $\omega$ (which we call primal weight) controls the balance between the primal and the dual variables. 

% We call $\omega\in(0,\infty)$ the primal weight and $s\in(0,\infty)$ the step-size. This allows us to control the scaling between the primal and dual iterates with a single parameter $\omega$ that weights the primal variables in the following norm
% \begin{equation*}
%     \|z\|_{\omega}:=\sqrt{w\|x\|_2^2+\frac{\|y\|_2^2}{\omega}}
% \end{equation*}
% This norm plays a role in the theory for PDHG and later algorithmic discussions.

\subsection{Algorithmic enhancements in PDLP}\label{sec:enhancement}
PDLP has essentially five major enhancements on top of restarted PDHG: presolving,  preconditioning, adaptive restart, adaptive step-size, and primal weight update.
\begin{itemize}
    \item {\bf Presolving.}
    PDLP utilizes PaPILO~\cite{gleixner2022papilo}, an open-sourced library, for the presolving step. The basic idea of presolving is to simplify the problem by detecting inconsistent bounds, removing empty rows and columns of the constraint matrix, removing variables whose lower and upper bounds are equal,  detecting duplicate rows and tightening bounds, etc.
    \item {\bf Preconditioning.}
    The performance of FOMs heavily depend on the condition number.
    PDLP utilizes a diagonal perconditioner to improve the condition number of the problem. More specifically, 
    it rescales the constraint matrix $K=(G,A)$ to $\tilde K=(\tilde G,\tilde A)=D_1KD_2$ with positive diagonal matrices $D_1$ and $D_2$, so that the resulting matrix $\tilde{K}$ is “well balanced”. Such preconditioning creates a new LP instance that replaces $A,G,c,b,h,u$ and $l$ in \eqref{eq:lp} with $\tilde G,\tilde A,\hat x=D_2^{-1}x,\tilde c=D_2c,(\tilde b,\tilde h)=D_1(b,h), \tilde u=D_2^{-1}u$ and $\tilde l=D_2^{-1}l$. For the default PDLP settings, a combination of Ruiz rescaling \cite{ruiz2001scaling} and the preconditioning technique proposed by Pock and Chambolle \cite{pock2011diagonal} is applied. 
    \item {\bf Adaptive restarts.} PDLP utilizes an adaptive restarting scheme that is similar (but not identical) to the one described below: essentially, we restart the algorithm whenever the normalized duality gap has a constant-factor shrinkage (we use $1/2$ below for simplicity):
\begin{equation}\label{eq:adaptive-restart}
        \rho_{\|\bar z^{n,k}-z^{n,0}\|_{2}}(\bar z^{n,k}) \leq \frac{1}{2} \rho_{\|z^{n,0}-z^{n-1,0}\|_{2}}(\bar z^{n,0})\ .\; 
\end{equation}
The normalized duality gap for LP can be computed with a linear time algorithm, thus one can efficiently implement this adaptive scheme. Adaptive restarts can speed up the convergence of the algorithm, in particular, for finding high-accuracy solutions.

%     PDLP utilizes an adaptive restarting scheme that is similar to the one described in Section \ref{sec:restart}. Adaptive restarts can speed up the convergence of the algorithm, in particular, for finding high-accuracy solutions.

%     An adaptive restart scheme is proposed in \cite{applegate2022faster}: essentially, we restart the algorithm whenever the normalized duality gap has a constant-factor shrinkage (we use $e^{-1}$ below for simplicity):
% \begin{equation}\label{eq:adaptive-restart}
%         \rho_{\|\bar z^{n,k}-z^{n,0}\|_{2}}(\bar z^{n,k}) \leq e^{-1} \rho_{\|z^{n,0}-z^{n-1,0}\|_{2}}(\bar z^{n,0})\ .\; 
% \end{equation}
% The normalized duality gap for LP can be computed with a linear time algorithm, thus one can efficiently implement this adaptive scheme.
    
%     In PDLP, we adaptively restart the PDHG algorithm in each outer iteration. As shown in Section \ref{sec:restart}, restart scheme enables PDHG achieving optimal rate for LP and it is numerically beneficial as well. The adaptive restart scheme discussed in Section \ref{sec:restart} is applied: restart happens when the normalized duality gap has sufficient decay, and more formally
%     \begin{equation*}
%     \begin{cases}
%         \rho_{\|\bar z^{n,k}-z^{n,0}\|_{2}}(\bar z^{n,k}) \leq e^{-1} \rho_{\|z^{n,0}-z^{n-1,0}\|_{2}}(\bar z^{n,0}),\; n\geq 1 \\ 
%         k \geq \tau^0,\; n=0
%     \end{cases} \ . 
% \end{equation*}

\item {\bf Adaptive step-size.} The theory-suggested step-size $1/\|A\|_2$ turns out to be too conservative in practice. PDLP tries to find a step-size by a heuristic line search that satisfies
\begin{equation}\label{eq:step-size}
        s\leq \frac{\|z^{k+1}-z^k\|_{\omega}^2}{2(y^{k+1}-y^k)^TK(x^{k+1}-x^k)}\ ,
    \end{equation}
where $\|z\|_{\omega}:=\sqrt{w\|x\|_2^2+\frac{\|y\|_2^2}{\omega}}$ and $w$ is the current primal weight. More details of the adaptive step-size rule can be found in \cite{applegate2021practical}. The inequality \eqref{eq:step-size} is inspired by the $\mathcal O(1/k)$ convergence rate proof of PDHG~\cite{chambolle2016ergodic,lu2023unified}. The adaptive step-size rule loses theoretical guarantee of the algorithm, but it reliably works in our numerical experiments.

    \item {\bf Primal weight update.}
    The primal weight $\omega$ aims to balance the primal space and the dual space with a heuristic fashion, and it is updated infrequently, only when restarting happens. The detailed description of primal weight update can be found in \cite{applegate2021practical}.
    
    % More specifically $\omega$ is initialized using 
    % \begin{equation*}
    % \mathrm{InitializePrimalWeight}(c,q):=\begin{cases}
    %     \frac{\|c\|_2}{\|q\|_2},\; & \text{if } \|c\|_2,\|q\|_2>\epsilon_{\mathrm{zero}}\\
    %     1, \; & \mathrm{otherwise}
    % \end{cases}
    % \end{equation*}
    % where $\epsilon_{\mathrm{zero}}$ is a small nonzero tolerance. Denote $\Delta_x^n=\|x^{n,0}-x^{n-1,0}\|_2$ and $\Delta_y^n=\|y^{n,0}-y^{n-1,0}\|_2$. PDLP applies primal weight update at the initiation of each new epoch.
    % {\small
    % \begin{equation*}
    % \mathrm{PrimalWeightUpdate}(z^{n,0},z^{n-1,0},\omega^{n-1}):=\begin{cases}
    %     \exp\pran{\theta \log\pran{\frac{\Delta_y^n}{\Delta_x^n}}+(1-\theta)\omega^{n-1}},\; & \Delta_x^n,\Delta_y^n>\epsilon_{\mathrm{zero}}\\
    %     \omega^{n-1}, \; & \mathrm{otherwise}
    % \end{cases}
    % \end{equation*}
    % }
    % The scheme aims to choose the primal weight $\omega^n$ such that distance to optimality in the primal and dual is the same, i.e., $\|(x^{n,k}-x^*,0)\|_{\omega^n}\approx \|(0,y^{n,k}-y^*)\|_{\omega^n}$. Moreover, PDLP performs a exponential smoothing with parameter $\theta \in[0,1]$ to dampen ocillations.

\end{itemize}

\subsection{Numerical performance of PDLP}
To illustrate the numerical performance of PDLP, we here present two sets of computational results on LP benchmark sets using the Julia implementation (these results were presented in \cite{applegate2021practical}). The experiments were performed on three datasets: 383 instances from the root-node relaxation of \texttt{MIPLIB 2017} collection~\cite{gleixner2021miplib} (we dub MIP Relaxations), 56 instances from Mittelmann's benchmark set~\cite{mittelmannbenchmark} (we dub LP Benchmark), and \texttt{Netlib} LP benchmark~\cite{netlib} (we dub Netlib). The progress metric we use is the relative KKT error, i.e., primal feasibility, dual feasibility and primal-dual gap, in the relative sense (see \cite{applegate2021practical} for a more formal definition). 

\begin{figure}[h!]
\begin{subfigure}{.33\textwidth}
\centering\includegraphics[width=1.0\linewidth]{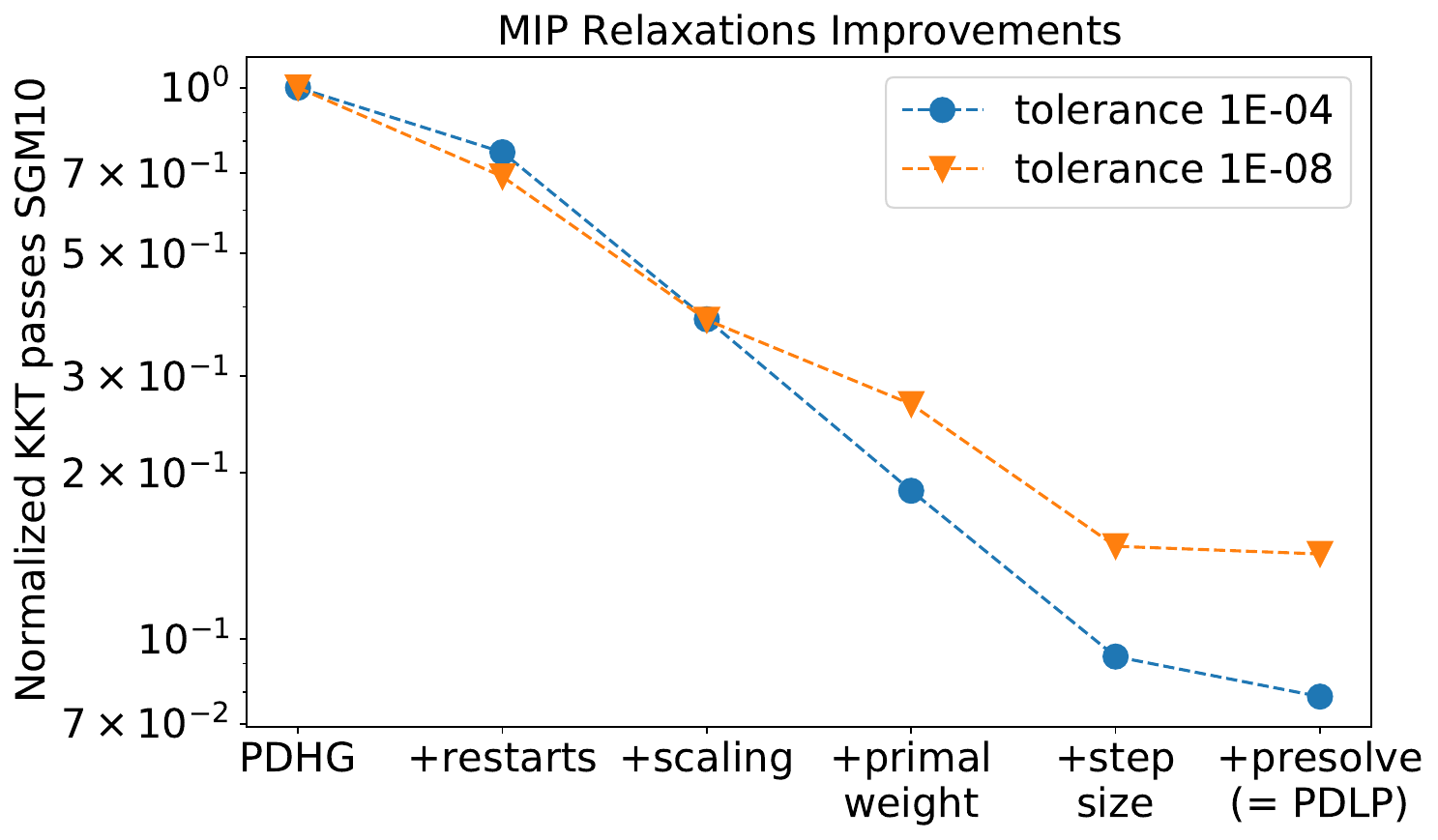}
\caption{\texttt{MIP Relaxations}}\label{fig:miplib-improvements}
\end{subfigure}
\begin{subfigure}{.33\textwidth}
\centering\includegraphics[width=1.0\linewidth]{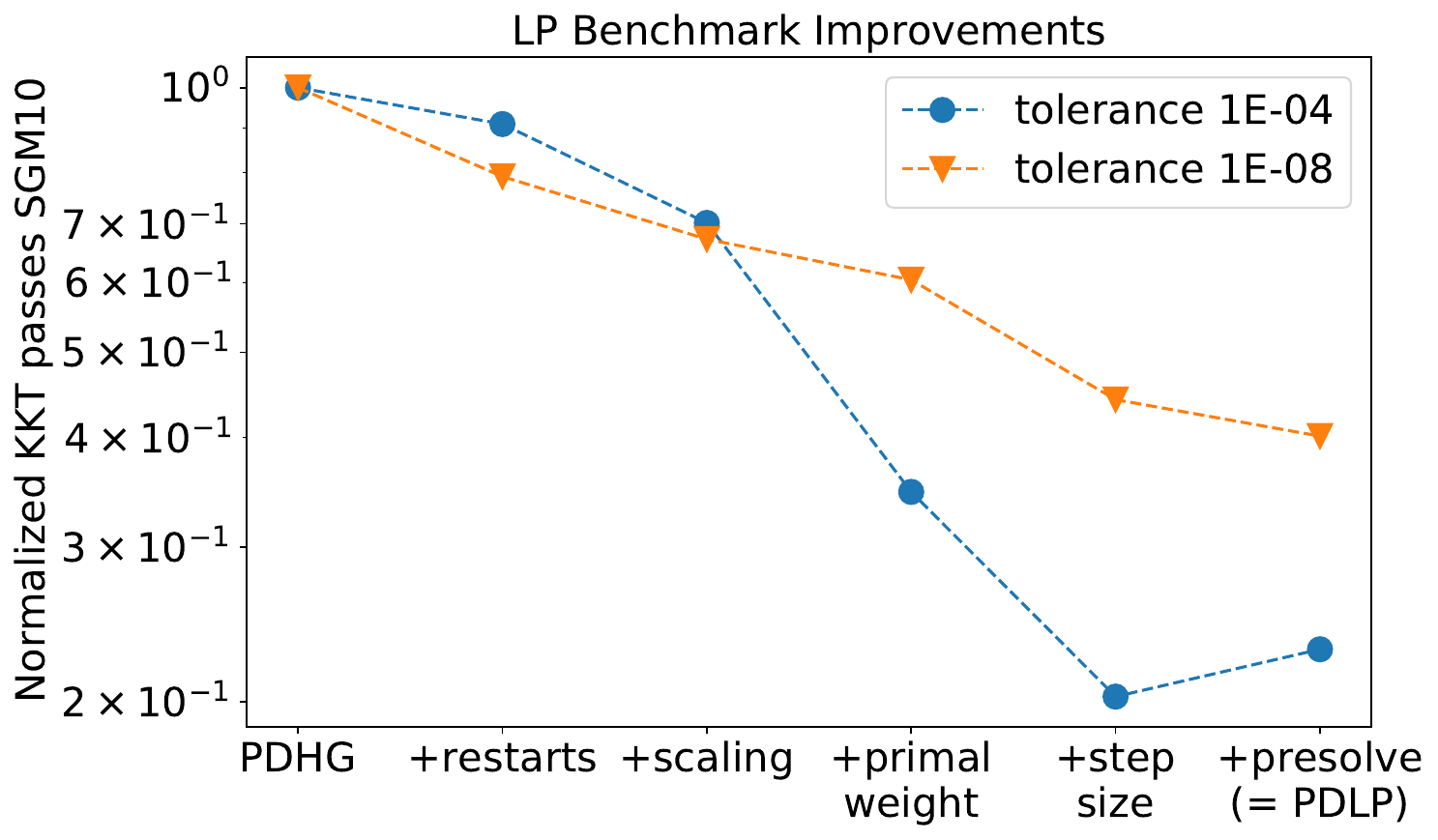}
\caption{\texttt{LP benchmark}}\label{fig:lpbenchmark-improvements}
\end{subfigure}
\begin{subfigure}{.33\textwidth}
\centering\includegraphics[width=1.0\linewidth]{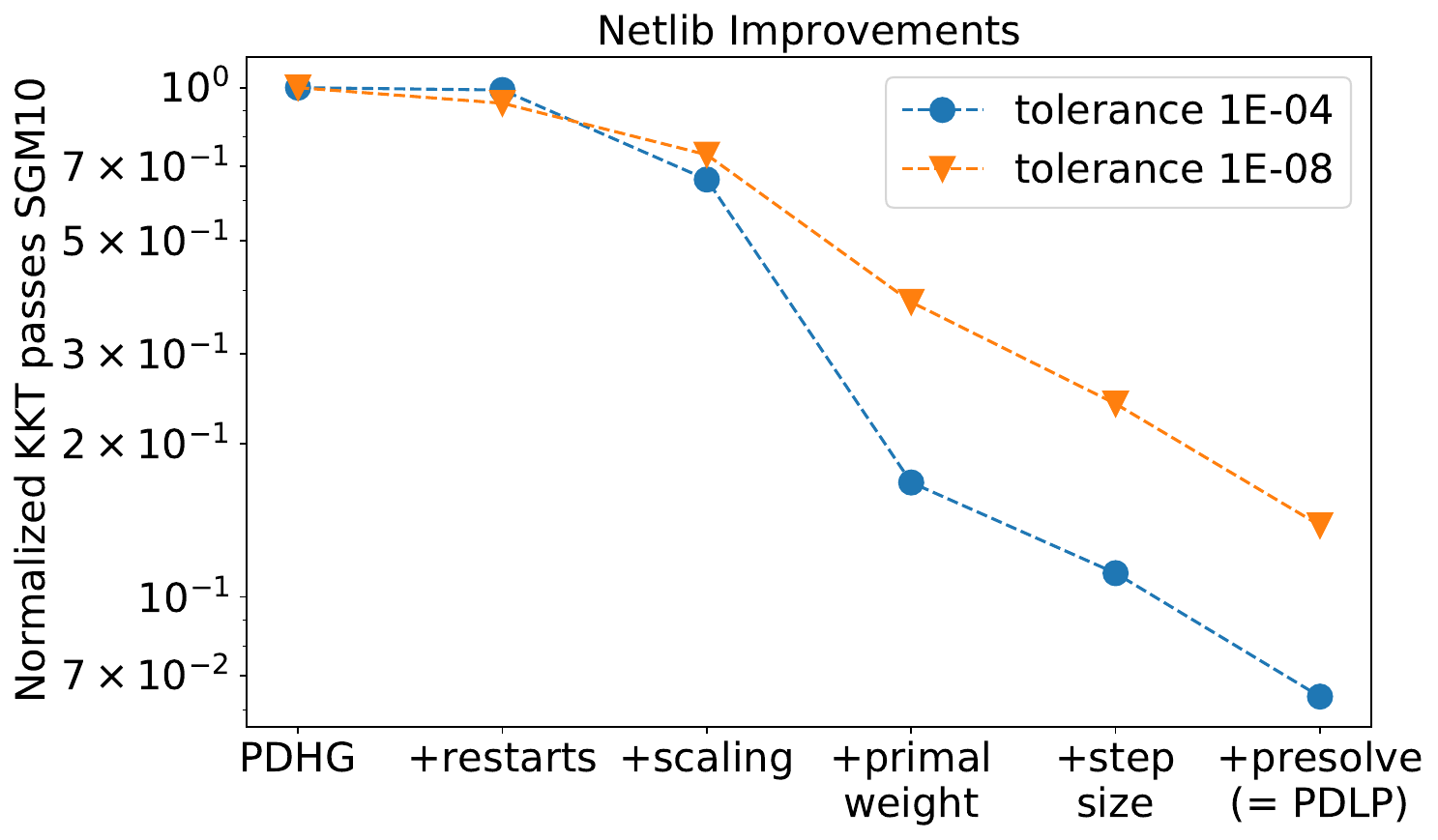}
\caption{\texttt{Netlib}}\label{fig:netlib-improvements}
\end{subfigure}
\caption{Summary of relative impact of PDLP's improvements}
\label{fig:pdlp-improvements}
\end{figure}

The first experiment is to demonstrate the effectiveness of the enhancements over the vanilla PDHG. Figure \ref{fig:pdlp-improvements} presents the relative improvements compared to vanilla PDHG by sequentially adding the enhancements. The y-axes of Figure \ref{fig:pdlp-improvements} display the shifted geometric mean (shifted by value 10) of the KKT passes normalized by the value for vanalla PDHG. We can see, with the exception of presolve for LP benchmark at tolerance $10^{-4}$, each of our enhancement in Section \ref{sec:enhancement} improves the performance of PDHG.

% 383 instances from \texttt{MIPLIB}, whose number of non-zeros are in between 100k to 10m. It looks at 

% We present two numerical experiments to study the effectiveness of PDLP with respect to traditional LP benchmark sets.

% First, we demonstrate PDLP's enhancements over baseline PDHG. The y-axis of Figure \ref{fig:pdlp-improvements} display the SGM10 of the KKT passes normalized by the value for baseline PDHG. We can see, with the exception of presolve for LP benchmark at tolerance $10^{-4}$, each of our enhancement in Section \ref{sec:enhancement} improves the performance of PDHG.

\begin{figure}[h!]
\begin{subfigure}{.33\textwidth}
  \centering
  \includegraphics[width=1.0\linewidth]{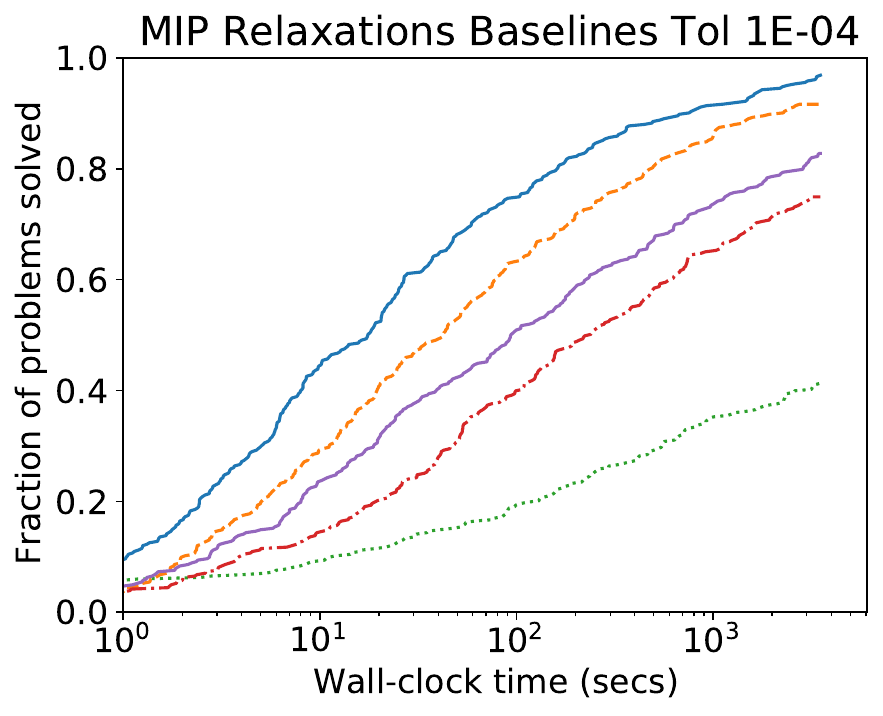}
%   \caption{1b}
\end{subfigure}
\begin{subfigure}{.33\textwidth}
  \centering
  \includegraphics[width=1.0\linewidth]{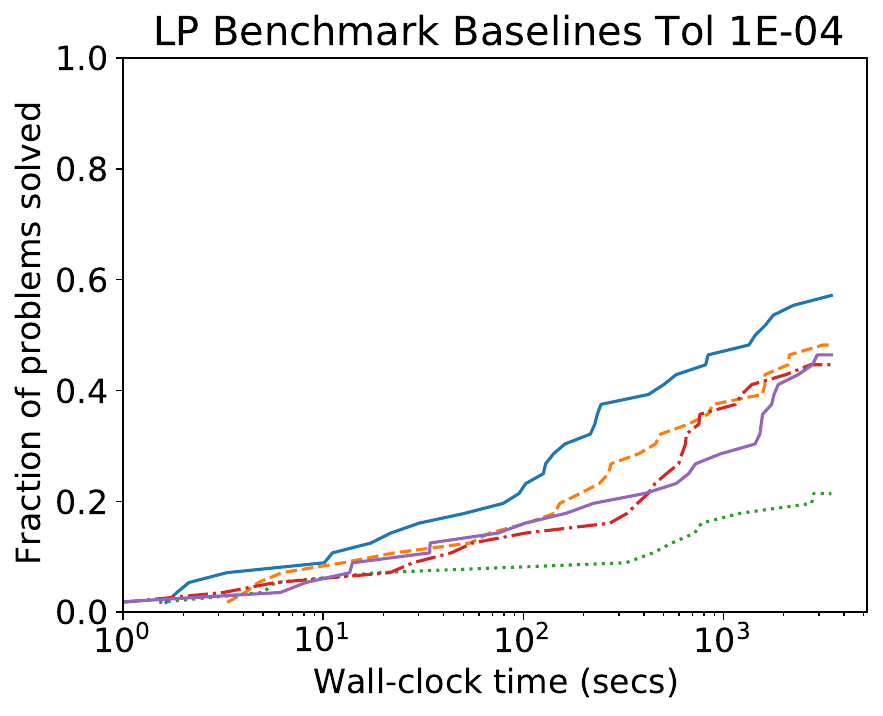}
%   \caption{1b}
\end{subfigure}
\begin{subfigure}{.33\textwidth}
  \centering
  \includegraphics[width=1.0\linewidth]{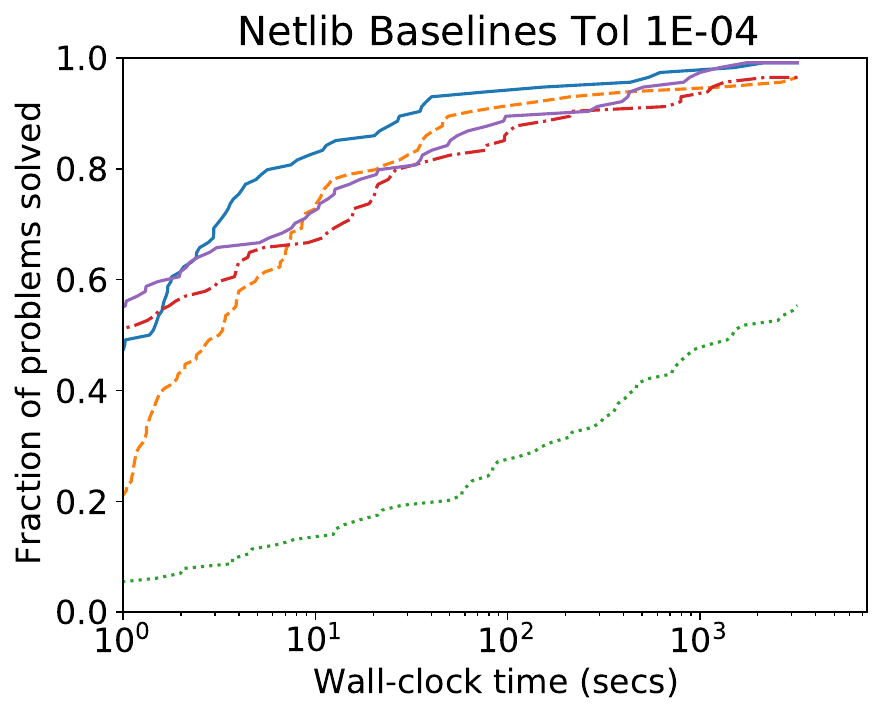}
%   \caption{1b}
\end{subfigure}\\
\begin{subfigure}{.33\textwidth}
  \centering
  \includegraphics[width=1.0\linewidth]{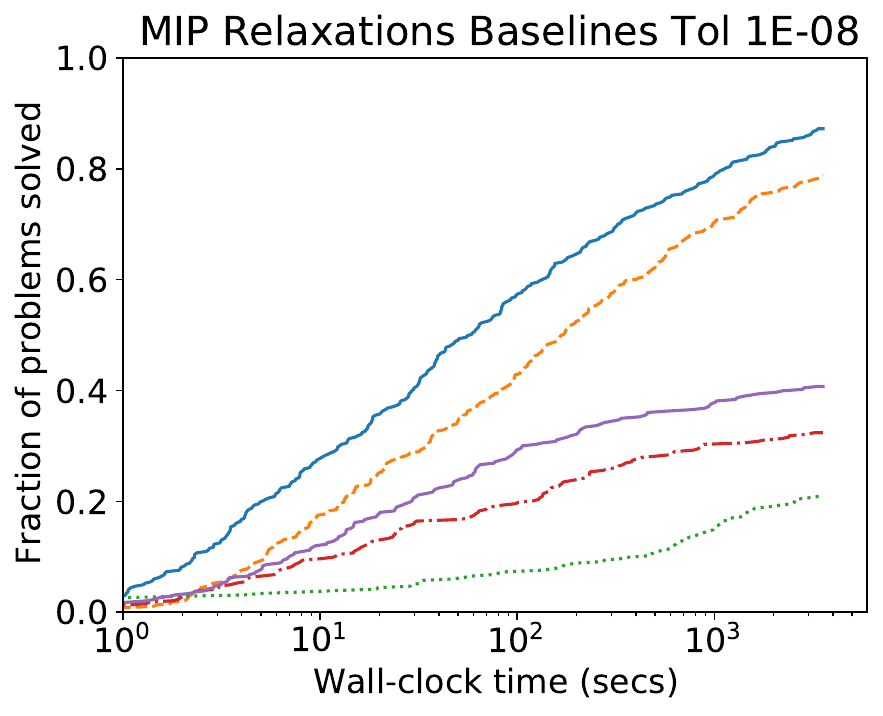}
\end{subfigure}
\begin{subfigure}{.33\textwidth}
  \centering
  \includegraphics[width=1.0\linewidth]{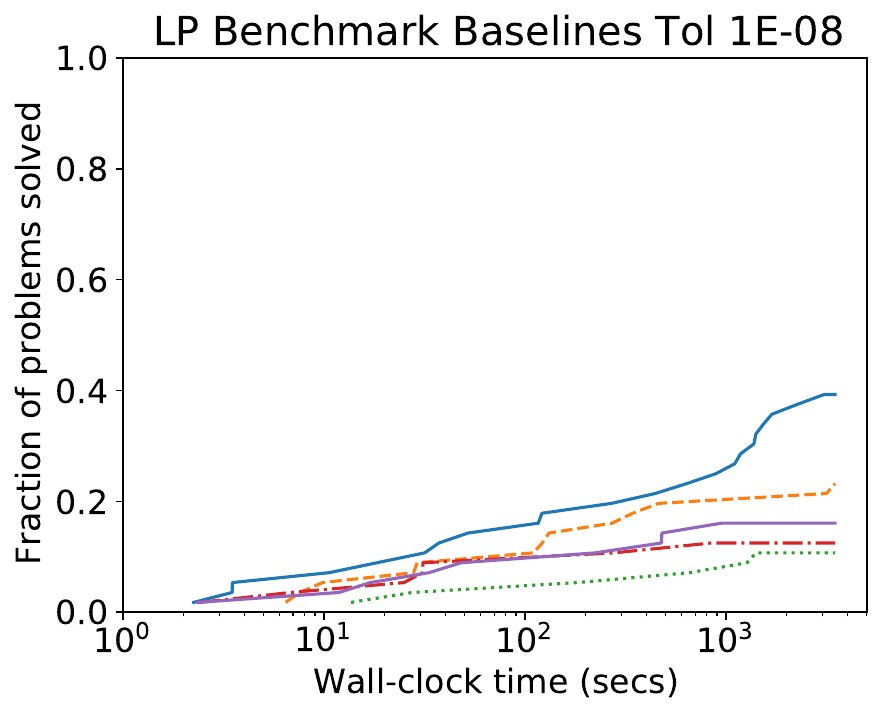}
  %\label{fig:sfig2}
\end{subfigure}
\begin{subfigure}{.33\textwidth}
  \centering
  \includegraphics[width=1.0\linewidth]{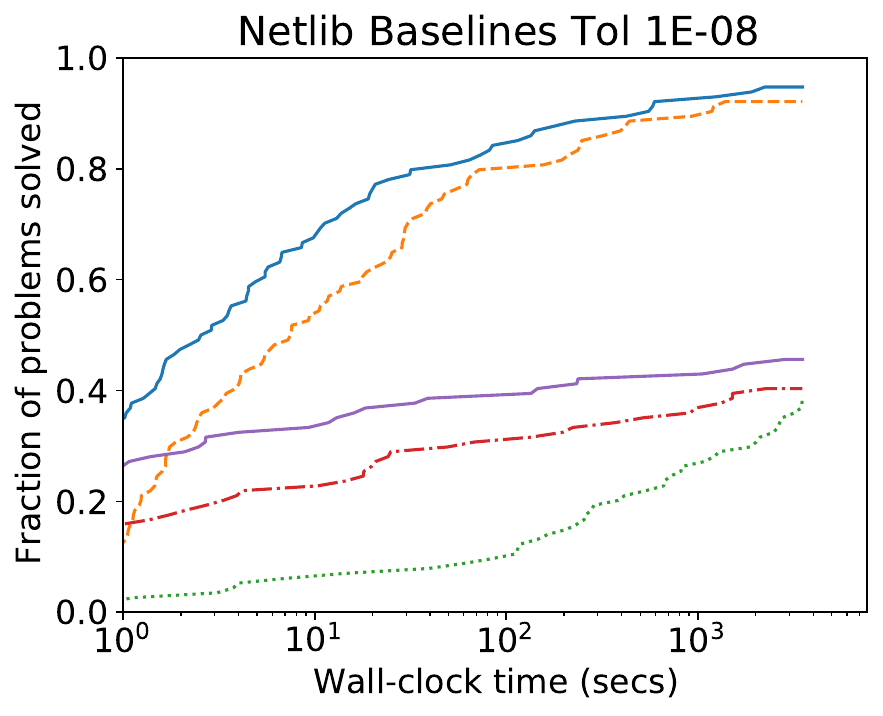}
  %\label{fig:sfig2}
\end{subfigure}\\
\begin{subfigure}{\textwidth}
\centering
\includegraphics[width=.95\linewidth]{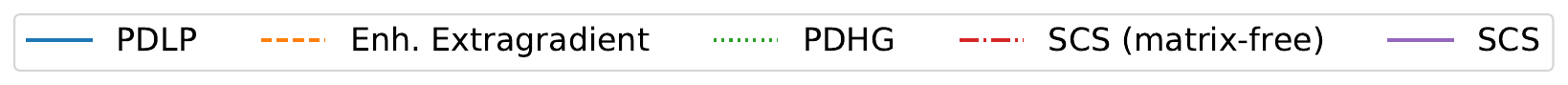}
\end{subfigure}
\caption{Number of problems solved for \texttt{MIP Relaxations} (left), \texttt{LP benchmark} (middle), and \texttt{Netlib} (right) datasets within one hour time limit.}
\label{fig:pdlp-vs-baseline}
\end{figure}
Figure \ref{fig:pdlp-vs-baseline} compares PDLP with other first-order methods:  SCS~\cite{o2016conic}, in both direct mode (i.e., solving the linear equation with factorization) and matrix-free mode (i.e., solving the linear equation with conjugate gradient), and our enhanced implementation of the extragradient method~\cite{korpelevich1976extragradient,nemirovski2004prox}. The extragradient method is a special case of mirror prox, and it is an approximation of PPM. Restarted extragradient has similar theoretical results as restarted PDHG~\cite{applegate2022faster}.
The comparisons are summarized in Figure \ref{fig:pdlp-vs-baseline}. We can see that PDLP have superior performance to both moderate accuracy $10^{-4}$ and high accuracy $10^{-8}$ compared with others.

\section{Case studies of large LP}
In this section, we present three case studies of PDLP on large LP instances (i.e., $\ge 10$m nonzeros) and compare its performance with Gurobi: personalized marketing, page rank, and robust production inventory problem. PDLP has superior performances in the first two instances, and Gurobi primal simplex has superior performance in the third instance. Overall, for the large instances that the factorization can fit in memory, PDLP won't completely substitute traditional LP solvers, it is reasonable to run them together in a portfolio. On the other hand, for problems that factorization cannot fit in memory, FOMs may be the only options.

% The first two instances

% we test the performance of PDLP on three large-scale applications of LP: PageRank, personalized marketing and robust production inventory problem. Note that LP formulations may not be the best approach to solve the problems but they are reasonable sources of very large instances.

\subsection{Personalized marketing}

Personalized marketing refers to companies sending out marketing treatments (such as discounts, coupons, etc) to individual customer periodically in order to attract more business. Such marketing treatments are usually limited by the total number of coupons that can be sent out, fairness considerations, etc. \cite{zhu2022augmented} proposes the problem using linear programming, and demonstrates the effectiveness of the model by utilizing the data of a departmental store collected from two states in the U.S..

More formally, suppose there are $I$ considered households and $J$ available marketing actions. The decision variables, $x_i^j\in [0,1]$, represent the probability a given household $i$ receives marketing action $j$. We denote the incremental profit that the firm earns from household $i$ if it receives marketing action $j$ as $r_i^j$ (these values can be learned from a machine learning model). The first set of constraints represents volume constraints on each marketing action, and $S_k$ is the set of households in customer segment $k$. The second set of constraints captures the volume constraint on all marketing actions. The combination of the marketing actions is determined by parameter $c_i^j$. The third set of constraints models the similarity constraints for each marketing action. The total number of households in customer segment $k$ is denoted by $n_k$, and $\lambda_j^{k_1k_2}$ restricts the difference between customer segments $k_1$ and $k_2$. The fourth set of constraints imposes similarity constraints on all marketing actions. The difference between customer segments $k_1$ and $k_2$ is restricted by $\gamma^{k_1k_2}$, and $d_i^j$ is the weighting factor to determine the combination of all marketing actions. The last two constraints restrict the firm’s action space so that each household has at most one marketing action. 
\begin{equation}
    \begin{aligned}
        \max_{x_i^j}& \ && \sum_{i=1}^I\sum_{j=1}^J r_i^jx_i^j\\
        \text{s.t.}& \  && a_k^j\leq \sum_{i\in S_k}x_i^j\leq b_k^j&&& \forall j\in [J], \;k\in [K]\\
        & \ && L_k\leq \sum_{i\in S_k}\sum_{j=1}^J c_i^jx_i^j\leq U_k&&& \forall k\in[K]\\
        & \ && \frac{1}{n_{k_1}}\sum_{i\in S_{k_1}}x_i^j\leq \lambda_j^{k_1k_2}\frac{1}{n_{k_2}}\sum_{i\in S_{k_2}}x_i^j &&& \forall j\in[J], k_1\in[K], k_2\in[K]\\
        & \ && \frac{1}{n_{k_1}}\sum_{i\in S_{k_1}}\sum_{j=1}^J d_i^jx_i^j\leq \gamma^{k_1k_2}\frac{1}{n_{k_2}}\sum_{i\in S_{k_2}}\sum_{j=1}^Jd_i^jx_i^j &&& \forall k_1\in[K], k_2\in[K]\\
        & \ && \sum_{j=1}^J x_i^j\leq 1 &&& \forall i\in[I]\\
        & \ && x_i^j\geq 0
    \end{aligned}
\end{equation}
Table~\ref{table:numberhousehold} summarizes the computational time of PDLP (C++ version) versus Gurobi (primal simplex, dual simplex, and barrier method) for four different models to $10^{-4}$ relative accuracy. For these instances, PDLP clearly shows its advantages, and it is often the case that Gurobi runs out of memory.
\begin{table}[]
			\centering
			\begin{threeparttable}
				%				\caption{Performance under Different Constraints.}\vspace{.2in}
				\begin{tabular}{ l c c c c c} 
					\toprule
					Model & \# nonzeros & \multicolumn{4}{c}{Computation Time of Gurobi / PDLP (in seconds)} \\    
					& & Gurobi Primal Simplex & Gurobi Dual Simplex  & Gurobi Barrier & PDLP \\
					%					& \multicolumn{4}{c}{Computation Time (in seconds)} \\
					\cline{3-6}
					%\hline
					%Volume A Upper& 15488 & 18469 & 16573 & 174\\
					model A& 25m & 15488 & 18469 & 16573 & 175\\
					model B & 37m & 31138 & - & - &  341\\
					model C & 25m & - & - & - & 136\\
					model D & 13m &- &- & - &  127 \\
%					\hline
%					&  \# nonzeros & \multicolumn{4}{c}{Optimal Profits (in 1,000)} \\
%					\cline{3-6}
%					%\hline
%					model A & 25m & 1121.02 & 1121.02 & 1121.02 & 1121.02 \\
%					model B & 37m & 1120.43 & - & - & 1120.43 \\
%					model C & 25m& - & - & - & 1202.87\\
%					model D & 13m & -&- & - &  751.26\\
					\bottomrule
				\end{tabular}
				\caption{Solve time in second of Gurobi and PDLP for four personalized marketing instances to $10^{-4}$ relative accuracy. ``-'' refers to raising an out-of-memory error.}\label{table:numberhousehold}
			\end{threeparttable}
		\end{table}

\subsection{PageRank}
\begin{table}[htbp]
\centering
\caption{Solve time for PageRank instances. Gurobi barrier has crossover disabled, 1 thread. PDLP and SCS solve to $10^{-8}$ relative accuracy. SCS is matrix-free. Baseline PDHG is unable to solve any instances. Presolve not applied. OOM = Out of Memory. The number of nonzero coefficients per instance is $8 \times (\text{\# nodes}) - 18$.}
\label{t:pagerank}
\begin{tabular}{lccccc}
\toprule
        \# nodes  &  PDLP & SCS & Gurobi Barrier  & Gurobi Primal Simp. &Gurobi Dual Simp.   \\
\midrule
$10^4$  &7.4 sec. & 1.3 sec. & 36 sec. & 37 sec. & 114 sec. \\
$10^5$  &35 sec. & 38 sec. & 7.8 hr. & 9.3 hr. & \textgreater24 hr.\\
$10^6$   &11 min. & 25 min. &OOM&  \textgreater 24 hr.& -\\
$10^7$  &5.4 hr. & 3.8 hr. &-&-& -\\
\bottomrule
\end{tabular}
\end{table}
PageRank refers to ranking web pages in the search engine results. There are multiple ways to model the problem, one of which is to  formulate the problem as finding a maximal right eigenvector of a stochastic matrix $S$ as a feasible solution of the LP problem~\cite{nesterov2014subgradient}, i.e.,
\begin{align}\label{eq:pagerank}
    \begin{split}
        & \mathrm{find}\;\; x \\
        & \ \mathrm{s.t.}\; Sx\leq x\\
        & \ \quad\;\; \;\mathbf{1}^Tx=1\\
        & \ \quad\;\; \;x\geq 0
    \end{split}
\end{align}
Nesterov~\cite{nesterov2014subgradient} states the constraint $\|x\|_{\infty}\geq 1$ to enforce $x\neq 0$. We instead use $\mathbf{1}^Tx = 1$ since it leads to a linear constraint.

\cite{applegate2021practical} generated a random scalable collection of pagerank instances, using Barab\'{a}si-Albert~\cite{barabasi1999emergence} preferential attachment graphs, and the Julia \texttt{LightGraphs.SimpleGraphs.barabasi\_albert} generator with degree set to 3. More specifically, an adjacency matrix is computed and scaled in the columns to make the matrix stochastic and this matrix is called
$S'$. Following the standard PageRank formulation, a damping factor to $S'$ is applied $S:=\lambda S'+(1-\lambda)J/n$ (where $J=\mathbf{1}\mathbf{1}^T$ is all-ones matrix). Intuitively, $S$ encodes a random walk that follows a link in the graph with probability $\lambda$ or jumps to a uniformly random node with probability $1-\lambda$. The direct approach to the damping factor results in a completely dense matrix. Instead, the fact that $Jx=1$ is used to rewrite the constraint $Sx\leq x$ in \eqref{eq:pagerank} as
\begin{equation*}
    \lambda(S'x)_i+(1-\lambda)/n\leq x_i,\;\forall i \ .
\end{equation*}

The results are summarized in Table~\ref{t:pagerank}. As we can see, when the instances get 10 times larger, the running time of PDLP and SCS roughly scale linearly as the size of the instance, whereas Gurobi barrier, primal simplex, and dual simplex scale much poorly. The fundamental reason for this is that the factorizations in the barrier and the simplex methods turn out to be much denser than the original constraint matrix. In this case, FOMs such as PDLP and SCS have clear advantages over simplex and barrier methods. We also would like to highlight that LP may not be the best way to solve PageRank problem, but that makes it easy to generate reasonable LPs of different sizes.

\subsection{Robust production inventory problem}

The production-inventory problem studies how to order products from factories in order to satisfy uncertain demand for a single product over a selling season. \cite{ben2004adjustable} proposes a linear decision rule for solving a robust version of this problem, and it initiates a new trend of research on robust optimization. The robust problem can be formulated as 

% One of the most popular approximation methods for solving multi-stage robust optimization problems is restricting the decisions in each stage to be a linear function of the uncertain variables observed in the past, as proposed in the seminal work of \cite{ben2004adjustable}. In the context of production inventory problem, the optimal linear decision rules can be obtained by solving the following linear programming~\cite{lu2022sparsity}

\begin{equation} \label{eq:ldr}
    \begin{aligned}
   & \underset{\substack{y_{t,1},\ldots,y_{t,t} \in \mathbb R^E:\;  \forall t \in [T]}}{\textnormal{min}} && \max_{\zeta_1 \in \mathcal{U}_1,\ldots,\zeta_{T+1} \in \mathcal{U}_{T+1}} \left \{ \sum_{t=1}^T \sum_{e=1}^E c_{te} \left( \sum_{s=1}^t    y_{t,s,e} \zeta_s\right) \right \} \\
    &\textnormal{subject to}&& \sum_{t=1}^T  \left( \sum_{s=1}^t y_{t,s,e} \zeta_s \right) \le  Q_e &&  \forall e \in[E] \\
   &&& 0 \le \left( \sum_{s=1}^t y_{t,s,e} \zeta_s \right) \le  p_{te} && \forall  e\in[E],\; t \in [T]\\
   &&&V_{\textnormal{min}} \le  v_1 + \sum_{\ell=1}^t \sum_{e=1}^E \left( \sum_{s=1}^\ell y_{\ell,s,e} \zeta_s \right)  - \sum_{s=2}^{t+1} \zeta_s  \le  V_{\textnormal{max}}&& \forall  t \in [T]\\
   &&& \quad \forall \zeta_1 \in \mathcal{U}_1,\ldots,\zeta_{T+1} \in \mathcal{U}_{T+1}\ ,
    \end{aligned}
\end{equation}
where $E$ is the number of factors, $T$ is the number of factories, $\xi_t$ is the uncertain demand that comes from a set $\mathcal{U}_t$, $y_{t,s,e}$ is the decision variable in the linear decision rule, $c_{te}$ is the per-unit cost, $Q_e$ is the  maximum total production level of factory $e$, $p_{te}$ is the maximum production level for factory $e$ in time period $t$, $V_{\min}$ and $V_{\max}$ specifies a remaining inventory level at each time period. This problem can be further formulated as linear programming by dualizing the inner maximization problem.

We compared the numerical performance of PDLP and Gurobi on a robust production inventory LP instance with $19$ millions non-zeros. While Gurobi dual simplex and barrier method failed to solve the problem within 1 days time limit, primal simplex solves it to optimality with 531160 iterations in 11800 seconds, and PDLP solves the problem to $10^{-4}$ accuracy with 358464 iterations in 26514 seconds. For this instance, Gurobi primal simplex outperforms PDLP. Our intuition is that this problem turns out to be highly degenerate, and there is a huge optimal solution space. Primal-simplex just needs to find one of the optimal extreme points, which can be efficient.

% The result is summarized as follows: Gurobi primal simplex solve to optimality with 531160 iterations in 11800 seconds, while PDLP solves to $10^{-4}$ accuracy with 358464 iterations in 26514 seconds. It turns out that PDLP can have inferior performance (more computing time) than Gurobi on some large-scale instances.

\section{Open questions}
It is an exciting time to see how FOMs can significantly scale up LP, an optimization problem that has been extensively studied since 1940s. This is indeed just the beginning time of FOMs for LP. There are still many open questions in this area, and we here mention three of them. 

First, while \cite{applegate2022faster} presents an optimal FOM for LP, which matches with the complexity lower bound, the linear convergence rate depends on Hoffman's constant of the KKT system, which is known to be exponentially loose (since it takes minimum of exponentially many items~\cite{pena2021new}), and clearly cannot characterize the numerical success of PDLP. A natural question is that, can we characterize the global convergence of PDHG for LP without using Hoffman's constant, or in other word, what is the fundamental geometric quantity of the LP instance that drives the convergence of PDHG (or more generally, FOMs)? 

Second, the state-of-the-art solvers for other continuous optimization problems, such as quadratic programming, second-order-cone programming, semi-definite programming, nonlinear programming, are mostly interior-point based algorithms. While SCS can also be used to solve some of these problems, it still needs to solve linear equations. How much can the success of FOMs for LP be extended to other optimization problems? What are the ``right'' FOMs for these problems? 

Third, how can FOMs be used to scale up mixed-integer programming (MIP)? In order to utilize the solution of FOMs for LP in the branch-and-bound tree, it is favorable to obtain an optimal basic feasible solution (BFS) to the LP, while FOMs usually do not directly output a BFS. Is there an efficient approach to obtain an optimal BFS from the optimal solution FOMs return? Furthermore, the LPs solved in the branch-and-bound tree are usually similar in nature. How can we take advantage of the warm start of LP solutions therein? 

% mixed integer linear programming utilizes branch-and-bound tree, and solves LP 

% \begin{itemize}
%     \item condition number
%     \item mixed-integer programming
% \end{itemize}

\bibliographystyle{amsplain}
\bibliography{ref-papers}
\end{document}